\documentclass{owrart}

\def\Alt{{\mathcal A}}
\def\C{{\bf C}}
\def\EE{{\mathcal E}}
\def\Pr{{\bf P}}
\def\Q{{\bf Q}}
\def\R{{\bf R}}
\def\XX{{\mathcal X}}
\def\X{X}
\def\Z{{\bf Z}}
\def\Fbar{{\,\overline{\!F}}}
\def\Qbar{{\,\overline{\!\Q\!}\,}}
\def\be{\begin{equation}}
\def\ee{\end{equation}}
\def\MW{Mordell--\mbox{\kern-.12em}Weil}
\def\NeS{N\'{e}ron--Severi}
\def\II{\mathop{\rm II}\nolimits}
\def\NS{\mathop{\rm NS}\nolimits}
\def\disc{\mathop{\rm disc}\nolimits}
\def\GL{\mathop{\rm GL}\nolimits}
\def\Gal{\mathop{\rm Gal}\nolimits}
\def\Ness{\mathop{N_{\rm ess}}\nolimits}
\def\Br{\mathop{\rm Br}\nolimits}

\begin{document}

%% ----------------------------------------------------------------
%% Example abstract.
%%
%% Please use the environment talk, it has three obligatory and one
%% optional argument. The syntax is:
%% -----------------------
%% \begin{talk}[coauthors]{Name of the speaker}{Title of the talk}{Author Sorting Index}
%%      .....
%% \end{talk}
%% -----------------------
%% The names of coauthors will appear in form of :
%% (joint work with ...)
%%
%% The Author Sorting Index is only used for sorting the list of all authors
%% by their last name. It should not contain special characters like
%% accents, german umlaute, etc. (just replace \"a by a, \'a by a, etc.)
%%
%% Please use the standard thebibliography environment to include
%% your references, and try to use labels for the bibitems, which
%% are uniquely assigned to you in order to avoid conflicts with other authors.
%% ---------------------------------------------------------------------

\begin{talk}{Noam D. Elkies}
{Three lectures on elliptic surfaces and curves of high rank}

\noindent
Over the past two years we have improved several of the (\MW)
rank records for elliptic curves over~$\Q$ and nonconstant
elliptic curves over~$\Q(t)$.  For example, we found the first example
of a curve $E/\Q$ with $28$ independent points $P_i \in E(\Q)$
(the previous record was $24$, by R.~Martin and W.~McMillen 2000),
and the first example of a curve over~$\Q$ with \MW\ group
$\cong (\Z/2\Z) \oplus \Z^{18}$ (the previous rank record for
a curve with a {$2$-torsion point} was $15$, by Dujella 2002).
In these lectures we give some of the background, theory,
and computational tools that led to these new records
and related applications.

\vspace*{2ex}

{\bf I} Context and overview: the theorems of Mordell(--Weil) and Mazur;
the rank problem; the approaches of N\'{e}ron--Shioda and Mestre; 
elliptic surfaces and N\'{e}ron specialization; fields other than~$\Q$.

{\bf II} Elliptic surfaces and K3 surfaces:
the \MW\ and \NeS\ groups;
K3 surfaces of high \NeS\ rank and their moduli;
an elliptic K3 surface over~$\Q$ of \MW\ rank~$17$.
Some other applications of K3 surfaces of high rank and their moduli.

{\bf III} Computational issues, techniques, and results:
slices of Niemeier lattices;
finding and transforming models of K3 surfaces of high rank;
searching for good specializations.
Summary of new rank records for elliptic curves.

\vspace*{4ex}

\noindent
{\bf I Context and overview.}
Mordell (1922) proved that the set $E(\Q)$ of rational points
of an elliptic curve $E/\Q$ has the structure of an abelian group,
and that this group is finitely generated.  That is,
$E(\Q) \cong T \oplus \Z^r$, where $T$\/ is a finite abelian group
(the \textit{torsion group} of~$E$\/) and $r$ is the
\textit{rank}\/ of~$E$.\footnote{
  Weil later (1928) generalized this from $E/\Q$ to $A/K$\/
  for an arbitrary abelian variety~$A$\/ over a number field~$K$,
  which is why the group $E(\Q)$ and rank~$r$ are often called the
  ``\MW\ group'' and ``\MW\ rank'' even in the case covered by
  Mordell's original result.}
This raises the basic structural question:

\vspace*{2ex}

\noindent
\textit{Which groups arise as $E(\Q)$ for some elliptic curve~$E/\Q$?}

\vspace*{2ex}

\noindent Equivalently,

\vspace*{2ex}

\noindent
\textit{Which ordered pairs $(T,r)$ arise as the torsion and rank
of some elliptic curve~$E/\Q$?}

\vspace*{2ex}

Mazur's celebrated torsion theorem \cite{Mazur} answers
the questions of which torsion groups arise: 
the cyclic groups of order~$N$\/ for $1 \leq N \leq 10$ and $N=12$,
and the groups $(\Z/2\Z) \oplus (\Z/2N\Z)$ for $1 \leq N \leq 4$.
These are exactly the fifteen groups~$T$\/ for which there is
a rational modular curve parametrizing elliptic curves~$E$\/
with an embedding of~$T$\/ into~$E(\Q)$.
It is thus almost immediate that each of these fifteen groups
arises infinitely often; the deep part of Mazur's theorem is
the proof that when the modular curve has positive genus
it has no rational points other than ``cusps''
(which parametrize certain degenerate elliptic curves).

This leaves the question of which values of~$r$ arise for each of
the fifteen possible~$T$.  At present this question is hopelessly
difficult.  It is not even known whether infinitely many~$r$ arise;
equivalently, whether $\limsup_{E/\Q}(r)$ is infinite for
some (all) of those~$T$.  As long as this question remains intractable,
we also ask for which $(T,r_0)$ can we prove that $\limsup(r) \geq r_0$;
this is in some sense a more demanding question than finding large
individual values of~$r$, in that proving $\limsup(r) \geq r_0$
requires infinitely many curves, not a single lucky guess.

Our main theme is the use of K3 surfaces of high rank and their moduli
to get new records for these questions (and also to obtain some
other applications of explicit parametrizations of K3 surfaces).
For the remainder of this first lecture we outline how the quest for
curves of large rank naturally leads to elliptic surfaces,
and illustrate two important earlier approaches to the problem.

Essentially the only technique known for proving lower bounds on
$\limsup(r)$ (at any rate the only technique known for $r_0>2$)
is finding \textit{parametrized families}, that is, infinite families of
elliptic curves $E$\/ together with generically independent points
$P_1,\ldots,P_{r_0}$.

Paradigmatic example: given $(x_i,y_i)$ $(i=1,2,3)$,
solve the simultaneous linear equations for $a_2,a_4,a_6$ that make
  $  y_i^2 = x_i^3 + a_2 x_i^2 + a_4 x_i + a_6  $
for each $i=1,2,3$.
This yields an elliptic curve~$E$\/ with $3$ rational points ($x_i,y_i$).
Exercise: they are generically independent (that is, independent
when $E$\/ is considered as an elliptic curve over the field
$\Q(x_1,y_1,x_2,y_2,x_3,y_3))$.  Hint: \textit{any}\/ quadruple
$(E,P_1,P_2,P_3)$ ($E$\/ some elliptic curve, each $P_i$ on~$E$\/)
arises this way for some $x_i,y_i$ if and only if each $P_i \neq 0$
and $P_i \neq \pm P_j$ for $i\neq j$.  [Moreover, we can make
$x_i,y_i$ unique by requiring $(x_1,x_2)=(0,1)$; then
$(x_3,y_1,y_2,y_3)$ gives a birational parametrization of the
``\hbox{$3$-rd} power $\EE^3$ of the universal elliptic curve''\kern-.2ex.]
By a specialization theorem of N\'{e}ron
(\cite{Neron}, see also \cite[Ch.11]{Serre:MW}),
later sharpened by Silverman (more on this below), 
there are infinitely many choices of $(x_i,y_i) \in \Q^6$
for which $P_1,P_2,P_3$ remain independent on the curve $E/\Q$.
Indeed (and not surprisingly),
this is true for ``most'' rational $(x_i,y_i)$,
and infinitely many non-isomorphic curves~$E$\/ arise this way.
Hence $\limsup(r) \geq 3$.

One quickly sees ways to improve this beyond rank~$3$;
for instance, use the extended Weierstrass form
$
y^2 + a_1 x y + a_3 = x^3 + a_2 x^2 + a_4 x + a_6
$
(or even the same with $a_0 x^3$ instead of $x^3$); or, pass
a cubic plane curve through $9$ ``random'' points of ${\bf P}^2$.
These are still not that far from each other, but they do suggest
two complementary ways of viewing the task.  In general, an elliptic
curve over $F(t_1,...,t_n)$ with several rational points is both:

1) [``\`{a} la Mestre''] Polynomial identities
(algebra, often ingeniously applied), and

2) [``\`{a} la N\'{e}ron''] An algebraic variety of dimension $n+1$
equipped with a suitable map to $n$-dimensional space over~$F$\/
(algebraic geometry).

We next give a short table of record ranks of
nonconstant elliptic curves\footnote{
  An elliptic curve over $\Q(t_1,\ldots,t_n)$ is ``nonconstant'' if it
  is not isomorphic over $\Q(t_1,\ldots,t_n)$ with a curve over~$\Q$.
  Such a curve of rank~$r$ yields nonconstant curves of rank at least~$r$
  over $\Q(t_1,\ldots,t_m)$ for each positive $m<n$ by specialization,
  so it is enough to consider $n=1$.
  }
over $\Q(t)$, all but the first and last of whose rows are taken from
\cite[Table~3]{RS} and represent the Mestre-style algebraic approach:

\vspace*{2ex}

\centerline{
\begin{tabular}{c|c}
  Rank $\geq$ & Author(s) and year \\ \hline
  8, 9, 10    & N\'{e}ron (1952) \\
    11, 12    & Mestre (1991) \\
      13      & Nagao (1994) \\
      14      & Mestre, Kihara (2000--1) \\
(15,16,)17,18 & NDE (2006--7)
\end{tabular}
}

\vspace*{2ex}

This leap from~$14$ to~$18$,\footnote{
  Ranks $15$ and $16$ are in parentheses because we proved the existence
  of such curves in~2006 but did not compute them explicitly.
  }
and similar improvements for curves with nontrivial torsion,
is also the key ingredient (via specialization) of
the new record ranks for individual curves over~$\Q$.
Curiously these improvements are achieved by returning to N\'{e}ron's
geometric viewpoint but applying it to elliptic surfaces at the next
level of complexity: elliptic K3 surfaces rather than rational
elliptic surfaces.  We shall say more about this in the second lecture.
First we interpolate some comments about N\'{e}ron's family,
an example of Mestre's identities, and remarks on
elliptic curves and surfaces defined over number fields other than~$\Q$.

Recall that we obtained rank $\geq 3$ by birationally parametrizing
\textit{all}\/ elliptic curves $E$\/ with three rational points $P_1,P_2,P_3$,
a.k.a.\ the ``\hbox{$3$-rd} power $\EE^3$
of the universal elliptic curve''; and observed that
some higher powers can be likewise parametrized using
other models of elliptic curves, notably the unique cubic curve
passing through $9$ general points $P_1,\ldots,P_9$ in the plane.
This already gives $\limsup(r) \geq 9$, because there are various ways
to get a $10$th point $P_0$ that can serve as an origin, such as the
$9$th \hbox{base-point} of the pencil of cubics through $P_1,...,P_8$.

For $r=10$, and thus for larger $r$, it is no longer possible to
completely parametrize $\EE^r$ --- ultimately because the modular
form $\Delta = q \prod_{n=1}^\infty (1-q^n)^{24}$
yields a holomorphic \hbox{$11$-form} on $\EE^{10}$!
Nevertheless, N\'{e}ron used the geometry of cubic curves to find
(in effect) rational curves on $\EE^{10}$ that give rise to elliptic curves
over $\Q(t)$ with $10$ generically independent rational points.
We describe this construction in some detail because some of the ideas
will recur in the more complicated setting of elliptic K3 surfaces.

We follow the exposition in \cite{Shioda:Neron}.  Start with
$P_1,...,P_8$, and thus also the ninth \hbox{base-point} $P_0$.
Then blowing up $\Pr^2$ at $P_0,P_1,\ldots,P_8$ gives
a birational isomorphism of $\Pr^2$ with the pencil of
cubics through these nine points.  Choose $P_0,...,P_8$ on a
\textit{cuspidal}\/ cubic, say $\Gamma: Y^2 Z = X^3$,
and choose the coordinate~$t$\/ on the pencil so that
$\Gamma$ is the preimage of $t = \infty$.
For generic such $P_1,\ldots,P_8$
this surface has no reducible fibers, and so rank~$8$.
Parametrize $\Gamma$ by $A(u)=(u:1:u^3)$, so that
$A(u_1),A(u_2),A(u_3)$ are collinear iff $u_1+u_2+u_3 = 0$; in particular,
the line through $A(u)$ and $A(-u/2)$ is tangent to~$\Gamma$ at $A(-u/2)$.
Let $D_1,D_2,D_3$ be these tangents for $u_1,u_2,u_3$.
Each $D_i$ meets a generic curve $E_t$ of the pencil at $2$ points
other than $P_i = A(u_i)$; so we get a ``double section'':
a section defined over a double cover of the \hbox{$t$-line}.
Moreover, each of these double covers is rational,
and $t=\infty$ is a branch point.  So any two of them give
a \hbox{degree-$4$} cover of $\Pr^1$ by a rational curve,
and N\'{e}ron shows that the two new points\footnote{
  Four new points are visible, but the two points in each double section
  sum to an element of the known \hbox{rank-$8$} group.
  }
are independent, giving rank~$10$ over~$\Q(t)$.  For each $t\in\Q$\/
(with finitely many exceptions where $E$\/ degenerates), we obtain
by specialization an elliptic curve $E_t/\Q$ with $10$ rational points.
Since $E/\Q(t)$ is nonconstant, N\'{e}ron's specialization theorem yields
infinitely many choices of $t$\/ where these points remain independent
and the curves $E_t$ are pairwise non-isomorphic.  (Silverman later used
the canonical height to construct, given a curve $E/\Q(t)$ and
independent rational points $P_1,\ldots,P_n$, an effective bound~$H$\/
such that the specialized points on $E_t$ remain independent for all~$t$
not of the form $t_0/t_1$ with $t_0,t_1 \in \Z \cap [-H,H]$, proving
that the set of exceptions to independence is finite and effectively
computable.  See again \cite[Ch.11]{Serre:MW}.)  Using all three $D_i$
yields rank~$11$ over the compositum of three rational double covers
of~$\Q(t)$, all branched at $t=\infty$.  This compositum is the
function field of an elliptic curve, usually of positive rank
and thus giving infinitely many examples of elliptic curves
over~$\Q$ with $11$ rational points.  A variation of N\'{e}ron's
specialization theorem, or of Silverman's refinement, then shows that 
these include infinitely many distinct curves of rank at least~$11$
over~$\Q$.

But this did not quite give a nonconstant elliptic curve of rank
$\geq 11$ over~$\Q(t)$.  Such a curve was first constructed in
\cite{Mestre:11}, as follows.  Suppose we have distinct 
$x_1,\ldots,x_{12} \in \Q$, polynomials $A_2, A_3 \in \Q(X)$ of degrees
at most~$2, 3$ respectively, and a monic polynomial $R(X)$
of degree~$4$ whose graph \hbox{$Y=R(X)$} intersects the plane
cubic curve $C: Y^3 + A_2(X) Y + A_3(X) = 0$ at the $12$ points
$P_i : (X,Y) = (x_i, R(x_i))$.  Then we expect to get rank $12-1$
by regarding $C$\/ as an elliptic curve with origin (say)~$P_1$.
Now the condition on the $x_i$, $A_j$, and~$R$\/ is equivalent to
$\prod_{i=1}^{12} (X-x_i) = R^3 + A_2 R + A_3$.  The $x_i$
thus uniquely determine $R$\/ as the principal part of the Taylor
expansion at infinity of $\left(\prod_{i=1}^{12} (X-x_i)\right)^{1/3}$,
and then we can recover $A_2$ and $A_3$ if and only if the
$X^{-1}$ coefficient of $\left(\prod_{i=1}^{12} (X-x_i)\right)^{1/3}$
vanishes (in which case $A_2,A_3$ are unique).  That coefficient
is a homogeneous quintic $F(x_1,\ldots,x_{12})$, which is
also invariant under translation $(x_i) \mapsto (x_i-\xi)$ and
thus yields \hbox{degree-$5$} hypersurface in $\Pr^{10}$.
This hypersurface contains some obvious rational subvarieties such as
the subspace cut out by $x_i + x_{6+i} = 0$ $(1 \leq i \leq 6)$,
but this choice makes our $11$ points dependent (though it gives $C$\/
a \hbox{$2$-torsion} point and can thus be used to construct
elliptic curves of moderately large rank with $T \supseteq \Z/2\Z$).
Mestre finds a less obvious rational subvariety of dimension~$3$ that
preserves independence, consisting of arbitrary linear combinations of
$(a,a,a,b,b,b,c,c,c,d,d,d)$ and $(b,c,d,a,c,d,a,b,d,a,b,c)$
for fixed $a,b,c,d$.\footnote{
  An extra dimension can be obtained by adding multiples of
  $(c,d,b,d,a,c,b,d,a,c,a,b)$ and $(d,b,c,c,d,a,d,a,b,b,c,a)$,
  the latter of which is redundant but highlights the $\Alt_4$ symmetry.
  This symmetry suggests the following equivalent construction
  of the resulting copy of $\Pr^2 \times \Pr^2$ in the hypersurface
  $F_5 = 0$: let $V$\/ be the irreducible \hbox{$3$-dimensional}
  representation of the alternating group~$\Alt_4$, let
  $\langle \cdot,\cdot \rangle$ be an \hbox{$\Alt_4$-invariant}
  perfect pairing on~$V\!$, and let $v,v'$ be any vectors in~$V$; then
  the $12$ inner products $\langle v, gv' \rangle$ $(g \in \Alt_4)$ 
  are coordinates $x_i$ of a point on~$F_5$.  This can be verified
  by regarding $F_5(x_1,\ldots,x_{12})$ as an \hbox{$\Alt_4^2$-invariant}
  polynomial on $V \oplus V$ and showing it must vanish identically
  [NDE 3.vii.1991, unpublished e-mail to J.-F.\ Mestre].
  }
Variations of this idea were later used by Mestre and others
[see table above] to push the record rank over $\Q(t)$ to~$14$
(and also for other purposes, such as constructing hyperelliptic curves
of given genus with many rational points).
Moreover, it was the \hbox{rank-$11$} family that was used by Mestre,
Nagao, Fermigier, Kouya, and Martin-McMillen during 1992--2000 to raise
the rank record for individual curves $E/\Q$ from $14$ to $15$, $17$,
$19$, $20$, $21$, $22$, $23$, and finally $24$ (see \cite{Dujella:hist}),
even after curves of rank $>11$ over $\Q(t)$ were found, because 
Mestre's curves have simpler coefficients and more parameters,
and thus offer greater scope for searching for high-rank
specializations.

What of the rank of elliptic curves $E/F(T)$ for general fields~$F$\/?
We exclude the case of a ``constant elliptic curve'', in which
$E$\/ is isomorphic over $F(T)$ with a fixed elliptic curve $E_0/F$,
because then $E(F(T))=E_0(F)$ (proof: no nonconstant maps from
$\Pr^1$ to $E_0$), which can be very large and/or complicated if
$F$\/ is large enough.  For nonconstant curves, the geometry of
the associated elliptic surface (more on this in the next lecture)
yields the result that $E(F(T))$ is finitely generated.  The list of
possible torsion groups is the same that Mazur proved over~$\Q$
together with $(\Z/N\Z)^2$ ($N=3,4,5$) and $(\Z/3\Z) \oplus (\Z/6\Z)$,
and the proof is much easier than over~$\Q$.  But,
as with Mordell's theorem for $E/\Q$, the bound on the rank is not uniform.
Indeed, when $F$\/ has characteristic $p>0$ the rank is unbounded
for ``isotrivial'' curves with constant and supersingular
\hbox{$j$-invariant} \cite{ST}, and also for non-isotrivial ones~\cite{Ulmer}.
In characteristic~zero, it is not yet known whether the rank
of nonconstant elliptic curves over $F(T)$ is unbounded even for $F=\C$;
the record is due to Shioda [1992]: the curve $y^2 = x^3 + T^n + 1$
over $\C(T)$ has (trivial torsion and) rank $\leq 68$, with equality
if and only if $360|n$.
Note that even though the curve is defined over~$\Q$,
most sections are not; for instance, if $3|n$ then $(x,y)=(-\mu T^{n/3},1)$
is on the curve for each $\mu\in\C$ such that $\mu^3=1$.
Still, the generators are all defined over some number field~$F_0$,
and it follows by N\'{e}ron's specialization theorem that there are
infinitely many elliptic curve of rank at least~$68$ over~$F_0$.

\vspace*{4ex}

\noindent
{\bf II Elliptic surfaces and K3 surfaces.}

[This part began with a review of the general setup of
elliptic curves~$E$\/ over $F(T)$ for an arbitrary field~$F \subset \C$,
relating the arithmetic of~$E$\/ with intersection theory on the
corresponding elliptic surface~$\XX$.  We do not repeat all of this
material here; see for instance \cite{Shioda:MW}.]

We assume throught that $\XX$\/ is a minimal N\'{e}ron model for~$E$.
Such a surface~$\XX$\/ is birational with an elliptic surface
in extended Weierstrass form
$
y^2 + a_1 x y + a_3 = x^3 + a_2 x^2 + a_4 x + a_6,
$
with each $a_i$ $(i=1,2,3,4,6)$ a section of $O(id)$
for some nonnegative integer~$d$\/ (in down-to-earth language,
a homogeneous polynomial of degree $id$\/ in two variables).
The smallest such~$d$\/ is the ``arithmetic genus'' of~$\XX$.
As the name suggests, the description of elliptic surfaces
of arithmetic genus~$d$\/ gets more complicated as~$d$\/ increases.
When $d=0$ we have a constant elliptic curve~$E_0$ over~$F(T)$
(equivalently, a surface~$\XX \cong E_0 \times \Pr^1$).
Once $d>0$, it follows from intersection theory on~$\XX$,
together with the fact that $h^{1,1}(\XX) = 10d$,
that the \MW\ rank of~$E/F(t)$ is at most $10d-2$.
Except for the smallest few~$d$\/ it is not known
whether this upper bound can be attained.

When $d=1$ we say $\XX$\/ is a ``rational elliptic surface''$\!$,
because it is birational with~$\Pr^2$,
at least over an algebraic closure $\Fbar$.
N\'{e}ron's surfaces of rank~$8$ are rational.
Since $8 = 10d-2$ for $d=1$, this gives the maximal \MW\ rank
of a rational elliptic surface.  Much more can be said of the
geometry and arithmetic of such surfaces, notably Shioda's
beautiful work relating rational elliptic surfaces with the
invariant theory of the Weyl group of the root lattice $E_8$
and its root sublattices; but we shall not follow this thread here.

Our main concern is the case $d=2$, when $\XX$\/ is an
``elliptic K3 surface''.  A K3 surface is a smooth algebraic surface
$\XX$\/ with trivial canonical class and $H^1(\XX,O_\XX)=0$.
This is the last case in which an algebraic surface can be elliptic
in more than one way; we heavily exploit this flexibility
in our analysis and computations.

A key invariant of a K3 surface~$\XX$\/ is its \textit{\NeS\ lattice}
$\NS(\XX) = \NS_\Fbar(\XX)$.  The \NeS\ lattice of any
compact algebraic surface over~$F$\/ is its \NeS\ group
(divisors defined over~$\Fbar$ modulo algebraic equivalence),
equipped with the symmetric integer-valued pairing induced from
the intersection pairing on divisors.  For a K3 surface,
this group is a free abelian group, and the pairing is \textit{even}:
$D \cdot D \in 2\Z$\/ for all $D \in \NS(\XX)$.
Let $\rho$ be the rank of~$\NS(\XX)$.  By the index theorem,
the pairing is nondegenerate of signature $(1,\rho-1)$.

If $\XX$\/ is elliptic then $\NS(\XX)$ contains two distinguished classes
defined over~$F$, the fiber~$f$\/ (preimage of any point under the map
$T: \XX \rightarrow \Pr^1$) and the zero-section~$s_0$.
The intersection pairing on the subgroup~$H$\/ they generate
is determined by \hbox{$f \cdot f = 0$}, $s_0 \cdot f = 1$, and
$s_0 \cdot s_0 = -d = -2$; hence $H$\/ is isomorphic with the
``hyperbolic plane'' (i.e., the even unimodular lattice with Gram matrix
{\small$\scriptscriptstyle{\left(\begin{matrix}0\!&\!1\cr1\!&\!0\end{matrix}\right)}$}).
Conversely, \textit{any} copy of~$H$\/ in $\NS(\XX)$ defined over~$F$\/
yields a model of~$\XX$\/ as an elliptic surface:
one of the generators or its negative is effective, and has
$2$ independent sections, whose ratio gives the desired map to $\Pr^1$.
(Warning: in general one might have to subtract some base locus
to recover the fiber class~$f$.)  Moreover, the pair $(\NS(\XX),H)$
determines the reducible fibers\footnote{
  Except for the distinctions between Kodaira types I$_1$ and~II
  (simple node and cusp, neither of which contributes to $\NS(\XX)$),
  I$_2$ and~III (either of which contributes $A_1$), and
  I$_3$ and~IV (either of which contributes $A_2$).
  }
and \MW\ group of the elliptic surface over~$\Fbar$.
Indeed let $\Ness$, the ``essential lattice'' of the elliptic surface,
be the orthogonal complement of~$H$\/ in $\NS(\XX)$,
with the pairing scaled by~$-1$ to make it positive definite.
Let $R \subseteq \Ness$ be the \textit{root lattice} of~$\Ness$,
the sublattice spanned by the roots (vectors of norm~$2$) in~$\Ness$.
This is a direct sum of root lattices
$A_n$ ($n\geq 1$), $D_n$ ($n\geq 4$), or $E_n$ ($n=6,7,8$),
with each factor indicating a reducible fiber of the corresponding type;
and the \MW\ group $E(\Fbar(T))$ is canonically isomorphic with
the quotient group $\Ness / R$.  In particular its rank is
the difference between the ranks of $\Ness$ and $R$.
The rank of $\Ness$, in turn, equals $\rho - 2$,
so the \MW\ rank is at most $\rho - 2$, with equality
if and only if $\Ness$ has no roots; in this case
the \MW\ rank over~$F(T)$ is also $\rho-2$ if and only if
$\NS(\XX)$ consists entirely of divisor classes defined over~$F$.

Now $\rho$ is at most $h^{1,1}(\XX) = 10d = 20$,
whence the upper bound $18 = 20 - 2$ on the \MW\ rank.
While a rational surface always has $\rho = h^{1,1}$,
for more complicated surfaces the \NeS\ rank~$\rho$
may be strictly smaller.  For K3 surfaces over~$\C$\/
the situation is completely described by the Torelli theorem
of Piateckii-Shapiro and \v{S}afarevi\v{c} [1971].
This theorem confirms and refines the following na\"{\i}ve
parameter count: there are $9+13-4 = 18$ parameters for
an elliptic K3 surface (the coefficients $a_4,a_6$ of
a narrow Weierstrass model have $8+1$ and \hbox{$12+1$} coefficients,
and we subtract $4$ for the dimension of $\GL_2$ acting on
the projective coordinates of the \hbox{$T$-line});
each free $\Ness$ generator, whether in~$R$\/ or the \MW\ group,
imposes one condition and thus reduces the dimension of
the moduli space by 1.  Recall that over~$\C$\/ it is known that
$H^2(\XX,\Z) \cong \II_{3,19}$ (the unique even unimodular lattice
of signature $(3,19)$), and that $\NS(\XX)$ embeds into $H^2(S,\Z)$.
The Torelli theorem asserts that this embedding is ``optimal'',
that is, realizes $\NS(\XX)$ as the intersection of $H^2(\XX,\Z)$
with a \hbox{$\Q$-vector subspace} of $H^2(\XX,\Z)\otimes\Q$;
for every such lattice~$L$\/ of signature $(1,\rho-1)$,
there is a nonempty (coarse) moduli space of pairs $(\XX, \iota)$,
where $\iota: L \rightarrow \NS(\XX)$ is an optimal embedding
consistent with the intersection pairing;
and each component of the moduli space has dimension $20-\rho$.
Moreover, for $\rho=20,19,18,17$ these moduli spaces repeat
some more familiar ones: CM elliptic curves for $\rho=20$,
elliptic and Shimura modular curves for $\rho=19$, and
moduli of abelian surfaces and RM abelian surfaces for certain cases of
$\rho=17$ and $\rho=18$.  It turns out that many of those moduli spaces
are more readily parametrized via K3 surfaces than by more direct approaches.
We shall treat these applications elsewhere, concentrating here on
the application to elliptic K3 surfaces.

To attain the upper bound of~$18$ on the \MW\ rank, we must use
a model of one of the (countably infinite number of) K3~surfaces
of \NeS\ rank~$20$ as an elliptic surface with trivial~$R$.
This can happen over~$\C$, and thus over~$\Qbar$ \cite{Cox,Nishiyama};
these proofs via \cite{PSS} use transcendental methods and yield
no explicit equations, but the example
\hbox{$Y^2 = X^3 - 27(T^{12} - 11T^6 - 1)$} was later obtained in \cite{CMT}.
This still leaves open the question of whether an elliptic K3 surface
can have \MW\ rank 18 over $\Q(T)$.  We repeat the warning that
it is not sufficient for the surface to be defined over~$\Q$;
as with Shioda's surface $Y^2 = X^3 + T^{360} + 1$, the
Chahal--Meijer--Top surface does not have all of $\NS(\XX)$
defined over~$\Q$.  Likewise for the \NeS\ groups
of some other familiar examples of K3 surfaces of maximal \NeS\ rank,
such as the diagonal quartic $X^4 + Y^4 = Z^4 + T^4$ in $\Pr^3$
or the complete intersection
$\sum_{i=1}^6 X_i = \sum_{i=1}^6 X_i^2 = \sum_{i=1}^6 X_i^3 = 0$
in $\Pr^5$.

One somewhat familiar example where the full \NeS\ group
\textit{is} defined over~$\Q$ is the universal elliptic curve with
a \hbox{$7$-torsion} point, considered naturally as an elliptic surface
over the modular curve $X_1(7) \cong \Pr^1_\Q$.  But this surface has
$\left|\disc(\NS(\XX))\right| = 7$, much too small for any of its
elliptic-surface models to have rank~$18$.  In fact we combine
arithmetic considerations with the construction in~\cite{Inose}
to show that if $\NS_\Q(\XX)$ has rank~$20$ then
$\disc(\NS(\XX))$ is one of the thirteen discriminants
$-3$, $-4$, $-7$, $-8$, $-11$, $-12$, $-16$, $-19$,
$-27$, $-28$, $-43$, $-67$, $-163$
of imaginary quadratic orders of class number~$1$.
Each of these arises uniquely up to twists, albeit with different
elliptic models --- already $-3$ and $-4$ have $6$ and $13$ respectively.
But even $163$ is too small for $\Ness$ to have no roots.\footnote{
  Such a lattice, positive-definite of rank~$18$ with discriminant~$163$
  and minimal norm $\geq 4$, would have broken the density record for
  a sphere packing in $\R^{18}$.  But the existence of such a lattice
  is not excluded by known sphere-packing bounds, so its impossibility
  had to be proved by other means.
  }
Therefore there are no elliptic K3 surfaces of \MW\ rank~$18$ over~$\Q$.\footnote{
  This was asserted in \cite{Shioda:K3}, but as a consequence of
  an incorrect result that was later retracted.
  }
But Mordell--Weil rank 17 is barely possible --- still not with
$\bigl(\rho,\, \left|\disc(\NS(\XX))\right|\bigr) = (20,163)$
but with an exceptional rational point on a certain Shimura curve!

More on this soon; first we describe torsion on elliptic K3 surfaces.
Each of the torsion groups in Mazur's list, other than
$\Z/9\Z$, $\Z/10\Z$, $\Z/12\Z$, and $(\Z/2\Z) \oplus (\Z/8\Z)$,
can arise for such a surface,
requiring at least the following reducible fibers,
and thus giving an upper bound on the rank equal to
$6$ less than the number of degenerate fibers
counted \textit{without} multiplicity:

\vspace*{1ex}

\centerline{
$
\begin{array}{c|ccccc}
\rm{torsion} &
\{0\} & \Z/2\Z & \Z/3\Z & \Z/4\Z & \Z/5\Z
\\ \hline
\rm{fibers} &
\phantom{1^{1^1}}
1^{24}
\phantom{1^{1^1}}
& 2^8 1^8 & 3^6 1^6 & 4^4 2^2 1^4 & 5^4 1^4 \\
\rm{formula} &
(0,0,0,a_4,a_6) & (0,a_2,0,a_4,0) & (a_1,0,a_3,0,0) &
  (a_1, a_2, a_1 a_2, 0, 0) & \rm {etc.} \\
\rm{bound} &
18 & 10 & 6 & 4 & 2
\end{array}
$
}

\vspace*{2ex}

\centerline{
$
\begin{array}{c|cccccc}
\rm{torsion} &
\Z/6\Z & \Z/7\Z & \Z/8\Z &
\! (\Z/2\Z)\oplus(\Z/2\Z) \! &
\! (\Z/2\Z)\oplus(\Z/4\Z) \! &
\! (\Z/2\Z)\oplus(\Z/6\Z) \!
\\ \hline
\rm{fibers} &
6^2 3^2 2^2 1^2 & 7^3 1^3 & 8^2 \, 4 \, 2 \, 1^2 &
\phantom{1^{1^1}}
2^{12}
\phantom{1^{1^1}}
& 4^4 2^4 & 6^3 2^3 \\
\rm{bound} &
2 & 0 & 0 & 6 & 2 & 0
\end{array}
$
}

\vspace*{2ex}

The three cases with bound zero are the universal elliptic curves
with that torsion group.  When the bound is positive it can always
be attained over~$\C$ but (as was already seen in the case of
trivial torsion) might not be attainable over~$\Q$.
The maximal rank is not known yet in each case,
because with nontrivial torsion it is possible for the
\MW\ group to be defined over~$\Q$ even though
$\Gal(\Qbar/\Q)$ acts nontrivially on $\NS_\Qbar(\XX)$.
Still, the discriminant $-163$ surface does have an elliptic model that
attains rank~$4$ with torsion group $\Z/4\Z$, and was used to get rank~$12$
over~$\Q$ (the previous rank record for an elliptic curve with a
rational \hbox{$4$-torsion} point was $9$, by Kulesz--Stahlke 2001).
Explicitly, the surface has equation $Y^2 + a X Y + a b Y = X^3 + a b X^2$
where $(a,b) = ((8T-1)(32T+7),$ $8(T+1)(15T-8)(31T-7))$;
it has a \hbox{$4$-torsion} point at $X=Y=0$, and four points with
$X = -15(T+1)(31T-7)(32T+7)/4$, $(8T-1)(15T-8)(31T-7)(32T+7)$,
$-(T+1)(8T-1)(15T-8)(32T+7)$, and $-4(T+1)(2T+5)(15T-8)(32T+7)$
that together with torsion generate $E(\Q(T))$; and taking
$T = 18745 / 6321$ yields a curve $E/\Q$\/ with eight further
independent points, so $E(\Q) \cong (\Z/4\Z) \oplus \Z^{12}$.
There are also various ways to combine pairs of quadratic sections to
get infinitely many $E/\Q$ with $E(\Q) \supseteq (\Z/4\Z) \oplus \Z^6$;
the previous rank record for an infinite family with \hbox{$4$-torsion}
was~$5$ (Kihara 2004, according to \cite{Dujella:generic},
where he cites two papers in \textit{Proc.\ Japan Acad.~A}).

A variant approach is to get some of the torsion group by
a suitable base change; for instance our records with torsion group
$(\Z/2\Z) \oplus (\Z/2\Z)$ were obtained by starting from an elliptic
K3 surface of \NeS\ rank~$20$ with torsion group $\Z/2\Z$
whose remaining \hbox{$2$-torsion} points are defined over
a quadratic extension of $\Q(T)$ that is still rational;
likewise we obtained $(\Z/2\Z) \oplus (\Z/4\Z)$
by quadratic base change from a curve over~$\Q(T)$
with torsion group $\Z/4\Z$.

We return now to the problem of finding elliptic K3 surfaces
of large \MW\ rank with no torsion restriction.
Having proved that rank~$18$ is unattainable,
we try for rank~$17$, corresponding to \NeS\ rank~$19$.
Here the moduli spaces have dimension \hbox{$20 - 19 = 1$,}
and in principle the Torelli theorem for K3 surfaces \cite{PSS}
identifies these curves with standard arithmetic quotients.
In practice it is not always easy to identify the modular curve
corresponding to a given lattice~$L$ of signature $(1,18)$,
especially when we need results over~$\Q$ rather than~$\C$.
But some identifications can be made.  For instance,
If $L \supset \II_{1,17}$ then the surfaces are parametrized by
the classical modular curve $\X_0(N)/w_N$ where $N = \disc(L)/2$.
This curve is rational for some rather large~$N$\/ (largest is $71$),
and elliptic of rank~$1$ for some~$N$\/ that are even larger
(largest is $131$).
For $N>131$ the curve has only finitely many rational points by
Mordell-Faltings.  We need only one rational point,
but it must be neither a cusp (because cusps yield degenerate surfaces)
nor a CM point (because CM points yield a surfaces of rank~$20$).
It is expected that there are only finitely many examples;
the largest known are for $N=191$ \cite{NDE:Atkin} and $N=311$
\cite{Galbraith}.  But again even those $N$\/ are too small
for $\Ness$ to have no roots.  Still, $N=311$ is large enough
for $R$\/ to have rank only~$2$, leaving \MW\ rank $17-2 = 15$.
This was already a new record, and as with N\'{e}ron's construction
it could be pushed a bit further with quadratic sections, to
$16$ over $\Q(T)$ and $17$ for infinitely many specializations.
(We can increment only once over $\Q(T)$, because for elliptic
K3 surfaces we cannot choose the ramification points.)
But I did not compute explicit equations for this K3 surface:
such a computation would have been a huge undertaking then,
and even now with better tools it would be a substantial project.
I did, however, manage to compute an elliptic model for the K3 surface
for the $N=191$ point that has a \hbox{$2$-torsion} point and
the minimal root lattice $A_1^8$ that can accommodate
$\Z/2\Z$ torsion.  Thus this elliptic surface has
\MW\ group $(\Z/2\Z) \oplus \Z^9$ over~$\Q(T)$.  Quadratic sections
increment this to $10$ over $\Q(T)$ and $11$ for an infinite family,
improving on Kihara's 2001 and 1997 records of~$9$
(\cite{Dujella:generic}, again citing papers in
\textit{Proc.\ Japan Acad.~A}).
Specialization of the K3 surface to $t\in\Q$ produced the new record
curve $E/\Q$ with $E(\Q) \cong (\Z/2\Z) \oplus \Z^{18}$.

One can do even better when $L$\/ is an even lattice of signature
$(1,18)$ that does not contain $\II_{1,17}$.  Let $N = \disc(L)/2$,
and suppose $N$\/ is squarefree.  Then $L \supset \II_{1,17}$
if and only if a certain obstruction in $\Br_2(\Q)$ vanishes.
This obstruction is supported on an even subset of the prime factors of~$N$.
If it does not vanish then we get the corresponding
Shimura modular curve instead of a classical (elliptic) modular curve.
When $N$\/ is composite, the Shimura curve can have smaller genus than
$\X_0(N)/w_N$ because there are fewer oldforms.  This lets us use
$N$\/ large enough that $\Ness$ can have trivial root system.
Even so, we did not find any case where the Shimura curve
has infinitely many rational points.  But for $N = 6\cdot 79$
we found a sporadic non-CM point.  Here the Shimura curve has
genus~$2$, and a bielliptic involution that lets us predict
an equation for the curve using the methods of~\cite{GonRot}.
We find $u^2 = 16 t^6 - 19 t^4 + 88 t^2 - 48$, with the following
rational points: the fixed points of the bielliptic involution
$(t,u) \leftrightarrow (-t,-u)$, with $t=\infty$;
four points with $|t|=2$ and $|u|=32$; and four with
$|t|=14/13$ and $|u| = 2^6 251 / 13^3$.  It turns out that the
last orbit is non-CM.  This one orbit of rational points yields
an elliptic K3 surface of \MW\ rank $17$ over~$\Q(t)$,
answering the question in~\cite{Shioda:K3} on the maximal
\MW\ rank of such a surface.  It also yields the new records of
$18$ for the \MW\ rank of a nonconstant elliptic curve over $\Q(T)$
(again via quadratic base change), and of $19$ for a lower bound on
$\limsup(r)$ over curves $E/\Q$.  Specialization of the \hbox{rank-$17$}
surface also yields several examples of elliptic curves over~$\Q$ with
more than $24$ independent rational points, including a curve
of rank at least~$28$.

Some remarks on the computation of these families and specializations
are in the third lecture.  We conclude this second lecture by noting
that the connection between K3 surfaces of \NeS\ rank~$19$ and
Shimura curves also makes it possible to compute explicit information
(equations, CM coordinates, Clebsch-Igusa invariants, etc.)\ about
Shimura curves of levels considerably beyond what was previously feasible.

\vspace*{4ex}

\noindent
{\bf III Computation and results.}

We briefly describe the steps of the computations needed to get from
the above theory of K3 surfaces and their moduli to explicit elliptic
curves over $\Q(T)$ and~$\Q$.

\textit{Finding suitable positive-definite lattices $\Ness$.}
After the second lecture's forced march through K3 territory,
I thought better of attempting another such review of Euclidean and
hyperbolic lattices.  Basically $\Ness$ is obtained as a suitable
slice of a Niemeier lattice.  The Niemeier lattices are the $24$
even unimodular lattices $\Lambda$ of rank $24$, each with a known
root system~$R$\/ and ``glue group'' $\Lambda/R$, which gives a handle
on the torsion and roots of its slices.
See \cite[Ch.~10.3, pages 399--402]{SPLAG}
for an example of this technique.

\textit{Finding equations for $E/\Q(T)$ and its \MW\ generators.}
Here it may seem that we are back where we started: we still seek
the coefficients of polynomial identities, such as
$y_i(T)^2 = x_i(T)^3 + a_4(T) x_i(T) + a_6(T)$ $(1 \leq i \leq 17)$,
with various auxiliary conditions on the $(x_i,y_i)$ to ensure
the correct height pairings.
There are too many variables to solve such a nonlinear system directly,
but in the \hbox{$4$-torsion} case shown earlier it was barely possible.
Still it was more convenient to eliminate only some of the variables,
and recover the remaining ones as follows.  The general theory
tells us that the coefficients are rational and behave well modulo
suitable small primes~$p$ such as $41 = (163+1)/4$.  An exhaustive
search mod~$p$ finds a solution.  Lift this solution arbitrarily
to characteristic zero and regard the lift as a \hbox{$p$-adic}
approximation to the correct solution.  Apply the natural generalization
of Newton's iteration $x \mapsto x - F(x)/F'(x)$ to this context,
using finite differences rather than derivatives to approximate~$F'$.
Each step doubles the \hbox{$p$-adic} precision.
Soon the \hbox{$p$-adic} approximation is close enough to recognize
the actual rational numbers by lattice reduction.  Then confirm them
by substitution into the desired identities.  Finally change coordinates
to simplify the equations to ones whose coefficients have smaller heights,
which is essential for finding high-rank specializations.

\textit{Exploiting different elliptic models of the same surface.}
Simple example: the Inose surfaces
$Y^2  =  X^3  +  A T^4 X  +  B'' T^7 + B T^6 + b' T^5$
over the \hbox{$T$-line} have essential lattice
$\Ness=R=E_8^2$ with reducible fibers at $T=0$ and $T = \infty$.
Scaling to \hbox{$Y^2  =  X^3  +  A X  +  B'' T + B + B'/T$}
we obtain an elliptic model over the \hbox{$X$-line}, this time with
$R = D_{16}$ and $[\Ness : R] = 2$
(note that $(T,Y)=(0,0)$ is a \hbox{$2$-torsion} point).
It turns out that the transformation is particularly simple
when, as here, the two lattices are \hbox{``$2$-neighbors''}:
they have isomorphic \hbox{index-$2$} sublattices.
We start from a model of~$\XX$\/ as an elliptic surface
whose coefficients are easier to compute, and then follow
a chain of \hbox{$2$-neighbors} (and the occasional \hbox{$3$-neighbor})
that introduces or removes roots and torsion until it reaches
an elliptic surface with the desired essential lattice.

\textit{Parametrizing families of K3 surfaces of \NeS\ rank $19$}.
When the \NeS\ rank is $19$ rather than $20$, our task is
not to solve for the coefficients of a single surface
but to parametrize a one-dimensional family by a modular curve.
We start at a known point $P_0$ of the curve
(maybe coming from a surface of rank $20$,
in an elliptic model in which $R$\/ is the same but the \MW\ rank
is larger by~$1$), and then deform it \hbox{$p$-adically}.
For example, fix a rational function $f$\/ of the coefficients
(say, the cross-ratio of the \hbox{$T$-coordinates} of
four reducible fibers), and use Newton to find points for which
$f$ is near $f(P_0)$.  Now the coefficients are generally
no longer rational even if $f(P)$ is, but they are algebraic with
degree at most $\deg(f)$.  We can guess those with lattice reduction
if $\deg(f)$ is small enough.  Varying $f(P)$ over
simple rational numbers in a $p$-adic neighborhood of~$f(P_0)$,
we can then guess the equations relating those coefficients with~$f$\/
by solving simultaneous linear equations.
At this point we have a guess for a (usually very singular)
model for the moduli curve, and if we can't or won't find a smooth model
directly then we can ask Magma (or someday Sage) for it.  We then verify
that the equations we guessed numerically actually work symbolically.
Then specialize to the \hbox{non-CM} point to find the desired surface.

\textit{Incrementing the rank via quadratic base change.}
As already noted, the necessary ``quadratic sections'' ---
rational curves on~$\XX$\/ that intersect the fiber~$f$\/ with
multiplicity~$2$ --- are harder to find than in N\'{e}ron's situation.
The trace of the quadratic section is an element of $E(\Q(T))$,
defined mod $2E(\Q(T))$ when we translate by elements of the \MW\ group;
equivalently, a half-lattice vector mod $E(\Q(T))$.
Intersection theory tells us that we need a coset of
$2E(\Q(T))$ in $E(\Q(T))$ consisting of vectors of norm~$2$ mod~$4$,
with no representatives of norm less than~$10$.
(This is for a surface with no reducible fibers; when
$R \neq \{0\}$ the criterion is more complicated.)
The corresponding coset of $E(\Q(T))$ in $\frac12 E(\Q(T))$
then consists of ``holes'' of norm at least $5/2$.
For our \hbox{rank-$17$} surface, there are literally thousands of
such holes, and for each one we get an inequivalent quadratic section.
The resulting \hbox{genus-zero} curves
are all rational because we can always find some other divisor whose
intersection with the quadratic section is odd.  This also gives
millions of biquadratic base changes of genus~$1$ and positive rank,
any one of which gets the lower bound $19$ on $\limsup(r)$.
(Alas none of them degenerates to a genus-zero curve, so we do not find
an elliptic curve of \MW\ rank at least~$19$ over~$\Q(T)$ this way.)

\textit{Guessing good specializations by Mestre's heuristic.}
The conjecture of Birch and Swinnerton-Dyer suggests that
large rank should correlate with small partial products of $L(E,1)$.
Taking logarithms, we want to make $\sum_{p<x} \log(p/N_p)$ very negative.
Experimentally, we need literally thousands of primes, and must canvass
many millions of specializations~$E_t$.  That's a lot of $N_p$'s to compute.
But each depends only on $t \bmod p$, so we precompute $\sum_{p<x} p$
of them once and for all, store a low-precision approximation to
$\log(p/N_p)$ for each one, and then search for large values of
$\sum_p \log(N_p/p)$ in sieve style.

\textit{Finding extra rational points.}  The resulting candidates
for large rank have coefficients much too large for it to be feasible
to find new rational points by direct search.  The simplest independent
set of $28$ rational points we could find on our record curve

\vspace*{2ex}

{
\small
  \centerline{
  $
  y^2 + xy + y = x^3 - x^2 -
  20067762415575526585033208209338542750930230312178956502\;x
  $
  }
  \vspace*{1ex}
  \centerline{\ \ 
  $
  {} + 34481611795030556467032985690390720374855944359319180361266008296291939448732243429 
  $
  }
}

\vspace*{2ex}

\noindent
has \hbox{$28$-digit} integers for its \hbox{$x$-coordinates}!
When the curve has nontrivial \hbox{$2$-torsion}, Cremona's
program \textsc{mwrank} quickly computes Selmer \hbox{$2$-groups}
to find upper bounds on the rank, and then usually finds enough
generators on the candidate record curves.  But in the absence
of torsion the coefficients are much too large for descent
to be feasible.  (This is why we can only say that the curve
has rank at least~$28$, not exactly~$28$, though it seems
quite unlikely that the rank is even larger.)  Instead we exploit
the known \hbox{rank-$17$} sublattice of $E/\Q$ to search for
rational points near half-lattice holes of the sublattice.
This yields equations $y^2 = Q(x)$ for quartics~$Q$\/ with
much smaller coefficients.  (This looks close enough to the behavior
of \hbox{$2$-descents} that the method might be regarded as a
fake \hbox{$2$-descent}.)  We then use a sieve technique,
implemented by C.~Stahlke and M.~Stoll in their C~program
\textsc{ratpoints}, to find a few such rational points
near some of the deepest half-lattice holes in the generic
\MW\ lattice.  Finally we use the canonical height to determine
the rank of the subgroup of $E(\Q)$ generated by all the known points.

\textit{Summary of new rank records.}
In the following table of record ranks of families of elliptic curves
with specified torsion group, ``$r+$'' means rank at least~$r$
over~$\Q(T)$, and at least $r+1$ for an infinite family parametrized by
a positive-rank elliptic curve obtained by quadratic base change from
the record curve over~$\Q(T)$.   Plain $r$ is a lower bound on the rank
of a curve over~$\Q(T)$.  The previous records as of 2004 are from
\cite{Dujella:generic}.

\vspace*{1ex}

\centerline{
\def\p{\phantom+}
$
\begin{array}{c|cccccc}
\rm{torsion} &
\{0\} & \Z/2\Z & \Z/3\Z & \Z/4\Z
   & \! (\Z/2\Z) \oplus (\Z/2\Z) \! & \! (\Z/2\Z) \oplus (\Z/4\Z) \!
\\ \hline
\leq 2004 &
14 & 9 & 6 & 5 & 6 & 3
\\
\rm{new} &
\p18+ & 10+ & 7 & \p5+ & \p7+ & \p3+
\end{array}
$
}

\vspace*{1ex}

For other torsion groups, the records remain
$3$ for $\Z/5\Z$ and $\Z/6\Z$;
``$1+$'' for
$\Z/7\Z$, $\Z/8\Z$, and $(\Z/2\Z) \oplus (\Z/6\Z)$
(the three cases where the universal elliptic curve is K3);
and ``$0+$'' for
$\Z/9\Z$, $\Z/10\Z$, $\Z/12\Z$, and $(\Z/2\Z) \oplus \Z/8\Z$
(the four cases where the universal elliptic curve has $d>2$).

The new rank records for individual curves over~$\Q$ are as follows:

\vspace*{1ex}

\centerline{
$
\begin{array}{c|ccccccc}
\rm{torsion} &
\{0\} & \Z/2\Z & \Z/4\Z & \Z/8\Z
   & \! (\Z/2\Z) \oplus (\Z/2\Z) \! & \! (\Z/2\Z) \oplus (\Z/4\Z) \!
   & \! (\Z/2\Z) \oplus (\Z/6\Z) \!
\\ \hline
\leq 2004 &
24 & 15 & 9 & 5 & 10 & 6 & 5
\\
\rm{new} &
28 & 18 & 12 & 6 & 14 & 8 & 6
\end{array}
$
}

\vspace*{1ex}

The incremental improvements for torsion groups
$\Z/8\Z$ and $(\Z/2\Z) \oplus (\Z/6\Z)$
are due only to better searching in known families.
The absence (so far?)\ of a new record for $\Z/3\Z$
may be due to the lack of an efficient implementation of
descent via a \hbox{$3$-isogeny}.

We conclude with a few remarks on integral points.
Our \hbox{$r \geq 28$} curve has at least $1174$ pairs $(x,\pm y)$
of integral points in its minimal model, but this is not a record:
a curve with $r \geq 25$ in the same family has at least $2810$ such pairs
in the known subgroup of $E(\Q)$.  The same family also contains
a curve for which we found only $21$ independent points but
the subgroup they generate contains at least $2564$ pairs 
of integral points.  Over $\Q(T)$, the analogue of integral points is
points $(X,Y)$ where $X,Y$\/ are polynomials of degree at most $4,6$
respectively.  In the absence of reducible fibers, these are exactly
the nonzero elements of the \MW\ group whose canonical height is
as small as~$4$.  Our elliptic K3 surface of rank~$17$ has
$1311$ such pairs $(X,\pm Y)$.

\end{talk}

\end{document}